\input amstex
\documentstyle{amsppt}
\magnification=\magstep1 \vsize=22 true cm \hsize=16.1 true cm
\voffset=1 true cm

\topmatter
\title    On angles in Teichm\"uller spaces
\endtitle
\title     On angles in Teichm\"uller spaces
\endtitle
\author   Yun Hu \&  Yuliang Shen
\endauthor
\affil (Department of Mathematics, Soochow   University)
\endaffil
\address Department of Mathematics,  Soochow  University,   Suzhou
215006,  P.\,R. China
\endaddress
\email\nofrills\eighttt huyun\_80\@163.com, ylshen\@suda.edu.cn
\endemail
\abstract We discuss the existence of the angle between two curves
in Teichm\"uller spaces and show that, in any infinite dimensional
Teichm\"uller space, there exist infinitely many geodesic triangles
each of which has the same three vertices and satisfies the property
that its three sides have the same and arbitrarily given length
while its three angles are equal to any given three possibly
different numbers from $0$ to $\pi$. This implies that the sum of
three angles of a geodesic triangle may be equal to any given number
from $0$ to $3\pi$ in an infinite dimensional Teichm\"uller space.
\endabstract
\thanks
\endthanks
\thanks
\endthanks
\thanks Research  supported by  the National Natural Science Foundation of China and the Specialized Research Fund for the Doctoral Program of Higher Education of China.
\endthanks
\endtopmatter

\document

\head 1 Introduction\endhead

The notion of geodesic segment plays an important role in the study
of the geometry of a metric space. Recall that a geodesic segment
 in a metric space is a continuous curve  such that for any
subarc  its length is equal to the distance between its two
endpoints. It is well known that there always exists a geodesic
segment between two points in any Teichm\"uller space (see [Ga]).
However, there are some essential differences of the geodesic
geometry between the finite and infinite dimensional Teichm\"uller
spaces (see [EKL], [EL], [Li1-4], [Sh], [Ta]). By Teichm\"uller's
theorem, there exists precisely one geodesic segment between two
points in finite dimensional Teichm\"uller spaces, while there exist
infinitely many geodesic segments joining certain pair of points in
any infinite dimensional Teichm\"uller space. The primary purpose of
the paper is to explore the further geodesic property of infinite
dimensional Teichm\"uller spaces.

It is  known that the inner product on the tangent space to a
Riemann manifold permits well-defined angle between two geodesic
segments. Since the Teichm\"uller distance is induced by a Finsler
structure (see [Ga], [Ob], [Ro]), it is not very clear how to define
the angle between two geodesic segments in a Teichm\"uller space.
Recently, following an idea of Professor Li, Yao [Ya] gave an
approach to define the angle between two geodesic segments in a
Teichm\"uller space, and showed that such an angle really exists in
a finite dimensional Teichm\"uller space. Later, Li-Qi [LQ] gave a
somewhat complicated condition under which there exists the angle
between two geodesic segments (of special form) in an infinite
dimensional Teichm\"uller space.

In this paper, we will continue to discuss the existence of the
angle between two geodesic segments in (infinite dimensional)
Teichm\"uller spaces. We first establish a variation formula for the
Teichm\"uller distance, from which it is proved that the angle
between two smooth geodesic segments exists in general. We then
study the geometry of Teichm\"uller spaces from the point of angle.
We show that in any infinite dimensional Teichm\"uller space, there
exist infinitely many geodesic triangles each of which has the same
three vertices and satisfies the property that  its three sides have
the same and arbitrarily given length while its three angles are
equal to any given three possibly different  numbers from  $0$ to
$\pi$. This implies that the sum of three
 angles of a geodesic triangle may be equal to any given number
from $0$ to $3\pi$ in an infinite dimensional Teichm\"uller space.
During the proof, we also find out that in any  infinite dimensional
Teichm\"uller space there do exist infinitely many  pairs of
geodesic segments between each pair of which the angle does not
exist. Consequently, in the view of angle, the geometry of an
infinite dimensional Teichm\"uller space is largely different from
the standard Euclidean or hyperbolic geometry.

\head 2 Preliminaries\endhead

In this section, we will recall some basic definitions and notations
from Teichm\"uller theory. For more details see the books [FM],
[Ga], [GL].

In what follows, $R$ will always denote  a hyperbolic Riemann
surface  covered by the unit disk in the complex plane. We denote
by $M(R)$ the unit ball of the space $L^{\infty}(R)$ of all
essentially bounded Beltrami differentials on $R$. We also denote by
$SQ(R)$ the unit sphere of the space $Q(R)$ of all integrable
holomorphic quadratic differentials on $R$.

For a given $\mu\in M(R)$, denote by $f^{\mu}$ the quasiconformal
mapping with domain $R$ and Beltrami coefficient $\mu$, which is
uniquely determined up to a conformal mapping on
$R^{\mu}=f^{\mu}(R)$. Two elements $\mu$ and $\nu$ in $M(R)$ are
equivalent, which is denoted by $\mu \sim\nu$, if $f^{\mu}$ and
$f^{\nu}$ are Teichm\"uller equivalent, meaning as usual that there
exists a conformal mapping $g$ from $R^{\mu}$ onto $R^{\nu}$ such
that $f^{\nu}$ and $g\circ f^{\mu}$ are homotopic (mod $\partial
R$). Then $T(R)=M(R)/\sim$ is the Teichm\"uller space of $R$. Let
$\Phi=\Phi_R: M(R)\rightarrow T(R)$ denote the canonical projection
from $M(R)$ to $T(R)$ so that $\Phi(\mu)$ is the equivalence
$[\mu]$. $\Phi(0)=[0]$ is called the base point of $T(R)$. It is
known that $T(R)$ is finite dimensional precisely when $R$ is of
finite type, namely, $R$ is a compact Riemann surface with possibly
finitely many points removed. It is also known that $T(R)$ has a
unique complex manifold structure such that $\Phi$ is a holomorphic
split submersion.

For any Beltrami coefficient $\mu\in M(R)$, define
$$k_0(\mu)=\inf\{\|\nu\|_{\infty}: \nu\sim\mu\},\tag 2.1$$
and set
$$K_0(\mu)=\frac{1+k_0(\mu)}{1-k_0(\mu)}.\tag 2.2$$
Then the Teichm\"uller distance $\tau=\tau_R$  between points
$\Phi(\mu_1)$ and $\Phi(\mu_2)$ is defined as
$$\tau (\Phi(\mu_1), \Phi(\mu_2))=\frac 12\log K_0(\mu),\tag 2.3$$
where $\mu$ is the Beltrami coefficient of the mapping
$f^{\mu_1}\circ (f^{\mu_2})^{-1}$. It is known that the
Teichm\"uller distance is compatible with the complex structure on
$T(R)$, namely, it is the Kobayashi metric on $T(R)$. We will need
an important fact about the Teichm\"uller distance: it is preserved
under a so-called allowable map. Recall that a Beltrami coefficient
$\mu$ in $M(R)$ induces an allowable map $A_{\mu}$ which maps $T(R)$
biholomorphically onto $T(R^{\mu})$ and sends $[\mu]$ to the base
point in $T(R^{\mu})$.

We say that $\mu\in M(R)$ is extremal if
$\|\mu\|_{\infty}=k_0(\mu)$. Then we also say that $f^{\mu}$ is
extremal. It is well known (see [Ha], [Kr], [RS] or Chapter 6 in
[Ga]) that $\mu$ is extremal if and only if $\mu$ satisfies the
Hamilton-Krushkal condition, that is, there exists a sequence
$(\phi_n)$ in $SQ(R)$ such that
$$\lim_{n\to\infty}\Re\iint_R\mu\phi_ndxdy=\|\mu\|_{\infty}.\tag 2.4$$
Such a sequence $(\phi_n)$ is called a Hamilton sequence for $\mu$. It is called degenerate if $\phi_n\to 0$ locally uniformly in $R$.

We also need a fundamental  inequality of Reich-Strebel (see [RS] or
Chapter 6 in [Ga]). We first introduce some notations. For any
$\mu\in M(R)$, set
$$I(\mu)=I_R(\mu)=\sup_{\phi\in
SQ(R)}\left|\Re\iint_R\frac{\mu\phi}{1-|\mu|^2}dxdy\right|,\tag
2.5$$
$$H(\mu)=H_R(\mu)=\sup_{\phi\in
SQ(R)}\left|\Re\iint_R\mu\phi dxdy\right|,\tag 2.6$$
$$J(\mu)=J_R(\mu)=\sup_{\phi\in
SQ(R)}\iint_R\frac{|\mu|^2|\phi|}{1-|\mu|^2}dxdy.\tag 2.7$$ Then, it
holds that
$$\frac{k_0(\mu)}{1-k_0(\mu)}-J(\mu)\le I(\mu)\le\frac{k_0(\mu)}{1+k_0(\mu)}+J(\mu).\tag
2.8$$

As stated in $\S 1$,  the Teichm\"uller distance is induced by a
Finsler structure (see [Ga], [Ob], [Ro]). For $\mu\in M(R)$ and
$\nu\in L^{\infty}(R)$, the Finsler structure $F=F_R$ is
$$F(\Phi(\mu),
\Phi'(\mu)\nu)=\inf\left\{\left\|\frac{\tilde\nu}{1-|\mu|^2}\right\|_{\infty}:\quad
\tilde\nu\in L^{\infty}(R) \text{ with $\Phi'(\mu)\tilde
\nu=\Phi'(\mu)\nu$}\right\}.\tag 2.9$$ From (2.8), it can be deduced
that
$$\align
F(\Phi(\mu),
\Phi'(\mu)\nu)&=H_{R^{\mu}}\left(\left(\frac{\nu}{1-|\mu|^2}\cdot\frac{{\partial}_z{f^{\mu}}}{\overline{{\partial}_z{f^{\mu}}}}\right)\circ(f^{\mu})^{-1}\right)\\
&=\sup_{\psi\in
SQ(R^{\mu})}\left|\Re\iint_{R^{\mu}}\psi\left[\left(\frac{\nu}{1-|\mu|^2}\cdot\frac{{\partial}_z{f^{\mu}}}{\overline{{\partial}_z{f^{\mu}}}}\right)\circ(f^{\mu})^{-1}\right]dudv\right|.\tag
2.10\endalign$$ In particular, $F(\Phi(0), \Phi'(0)\nu)=H(\nu)$. It
is known that the Finsler structure $F$ is continuous on the tangent
bundle of the Teichm\"uller space $T(R)$.

\head 3 A variation formula \endhead

 Let $\mu(t)$ be a continuous curve from $[0, t_0]$ into $M(R)$. We
 say $\mu(t)$ is differentiable at $0$ if there exist some $\mu\in
 L^{\infty}(R)$ such that $\mu(t)=\mu(0)+t\mu+o(t)$ as $t\to 0_+$,
 or more precisely,
 $$\lim_{t\to 0_+}\left\|\frac{\mu(t)-\mu(0)}{t}-\mu\right\|_{\infty}=0.\tag 3.1$$
We call $\mu$ the derivative of $\mu(t)$ at $0$, and denote it by
$\mu'(0)$.

\proclaim{Theorem 3.1} Let $\mu(t)$ and $\nu(t)$ be two continuous
curves from $[0, t_0]$ into $M(R)$ which are differentiable at $0$
and satisfy $\mu(0)=\nu(0)$. Then it holds that
$$\tau(\Phi(\mu(t)), \Phi(\nu(t)))=tF(\Phi(\mu(0)),
\Phi'(\mu(0))(\mu'(0)-\nu'(0)))+o(t),\quad t\to 0_+.\tag 3.2$$
\endproclaim

The following corollary is an immediate consequence of Theorem 3.1.
When $R$ is a compact Riemann surface, it was proved by Yao [Ya] by
a lengthy computation.

\proclaim{Corollary 3.1} For any two Beltrami differentials  $\mu$
and $\nu$ in $L^{\infty}(R)$, it holds that
$$\tau(\Phi(t\mu), \Phi(t\nu))=tH(\mu-\nu)+o(t),\quad t\to 0_+.$$
\endproclaim

To prove Theorem 3.1, we need the following lemma, which is a direct
consequence of the fundamental inequality (2.8).

 \proclaim{Lemma 3.1}
Suppose $R_t$ is a Riemann surface which may depend on $t\in [0,
t_0]$. If $\eta(t)\in M(R_t)$ satisfies $\eta(t)=t\delta(t)+o(t)$ as
$t\to 0_+$, where $\delta(t)\in L^{\infty}(R_t)$ satisfies
$\delta(t)=O(1)$ as $t\to 0_+$, then it holds that
$$\tau_{R_t}(\Phi_{R_t}(0), \Phi_{R_t}(\eta(t)))=tH_{R_t}(\delta(t))+o(t),\quad t\to 0_+.\tag 3.3$$
\endproclaim

\demo{Proof} By definition, $\tau_{R_t}(\Phi_{R_t}(0),
\Phi_{R_t}(\eta(t)))=(1+o(1))k_0(\eta(t))$ as $t\to 0_+$. Now we
replace $\mu$ by $\eta(t)$ in the inequality (2.8) on the Riemann
surface $R_t$. Clearly, $I_{R_t}(\eta(t))$ differs from
$H_{R_t}(\eta(t))$ by a term of order $t^2$, both
$k_0(\eta(t))/(1-k_0(\eta(t)))$ and $k_0(\eta(t))/(1+k_0(\eta(t)))$ differ
from $k_0(\eta(t))$ also by a term of order $t^2$, while
$J_{R_t}(\eta(t))$ is a term of ordered $t^2$. We conclude that
$k_0(\eta(t))=H_{R_t}(\eta(t))+o(t)$ as $t\to 0_+$ and (3.3)
follows.

\enddemo

\noindent {\it Proof of Theorem 3.1}\quad Let $\eta(t)$ be the Beltrami coefficient of
$f^{\mu(t)}\circ(f^{\nu(t)})^{-1}$, namely,
$$\eta(t)=\left(\frac{\mu(t)-\nu(t)}{1-\overline{\nu(t)}\mu(t)}\cdot\frac{{\partial}_z
f^{\nu(t)}}{\overline{{\partial}_z f^{\nu(t)}}}\right)\circ
(f^{\nu(t)})^{-1}.$$ By the differentiability of $\mu(t)$ and
$\nu(t)$, we obtain
$$\frac{\mu(t)-\nu(t)}{1-\overline{\nu(t)}\mu(t)}=\frac{t(\mu'(0)-\nu'(0))}{1-|\nu(t)|^2}+o(t),\quad
t\to 0_+.$$ Clearly, with $R_t=R^{\nu(t)}$, $\eta(t)\in M(R_t)$
satisfies the assumption of Lemma 3.1 with
$$\delta(t)=\left(\frac{(\mu'(0)-\nu'(0))}{1-|\nu(t)|^2}\cdot\frac{{\partial}_z
f^{\nu(t)}}{\overline{{\partial}_z f^{\nu(t)}}}\right)\circ
(f^{\nu(t)})^{-1}.$$ By Lemma 3.1, $$\tau_{R_t}(\Phi_{R_t}(0),
\Phi_{R_t}(\eta(t)))=tH_{R_t}(\delta(t))+o(t),\quad t\to 0_+.$$ But
by (2.10) and the continuity of the Finsler structure $F$, as $t\to
0_+$ it holds that
$$H_{R_t}(\delta(t))=F(\Phi(\nu(t)),
\Phi'(\nu(t))(\mu'(0)-\nu'(0)))\to F(\Phi(\mu(0)),
\Phi'(\mu(0))(\mu'(0)-\nu'(0))).$$ Thus,
$$\tau_{R_t}(\Phi_{R_t}(0), \Phi_{R_t}(\eta(t)))=tF(\Phi(\mu(0)),
\Phi'(\mu(0))(\mu'(0)-\nu'(0)))+o(t),\quad t\to
0_+.$$ Finally, $\tau(\Phi(\mu(t)),
\Phi(\nu(t)))=\tau_{R_t}(\Phi_{R_t}(0), \Phi_{R_t}(\eta(t)))$ by the
distance-preserving property of the allowable map $A_{\nu(t)}:
T(R)\to T(R_t)$ and (3.4) follows.

\head{4 Existence of angle}\endhead

We first introduce the notion of the angle between two joint curves
in a general metric space $(X, d)$. Let $\alpha$ and $\beta$ be two
continuous curves in $X$ with one  common endpoint $p$. For any
$r>0$, we choose $x(r)\in\alpha$ and $y(r)\in\beta$ such that the
length of the sub-curve of $\alpha$ between $p$ and $x(r)$ is the
same as that of the sub-curve of $\beta$ between $p$ and $y(r)$ and
equal to $r$. Then the angle at $p$ between $\alpha$ and $\beta$,
denoted by $\langle\alpha, \beta\rangle_p$, is defined as the number
$\theta\in [0, \pi]$ by the equation
$$2\sin\frac{\theta}{2}=\lim_{r\to 0}\frac{d(x(r),
y(r))}{r},\tag 4.1$$ if the limit exists. Notice that when both
$\alpha$ and $\beta$ are geodesic segments in a Teichm\"uller space,
the notion of the angle is reduced to the one  introduced by Yao
[Ya] and Li-Qi [LQ]. A trivial case is when $\alpha\equiv\beta$ in a
neighborhood of $p$, then the angle at $p$ between $\alpha$ and
$\beta$ exists and equals $0$. Another trivial case is when
$\alpha\cup\beta$ is geodesic at $p$, namely, there exists some
closed neighborhood $U(p)$ of $p$ such that $(\alpha\cup\beta)\cap
U(p)$ is a geodesic segment. Then, it is clear that the angle at $p$
between $\alpha$ and $\beta$ exists and equals $\pi$. In what
follows we always assume that $\alpha\neq\beta$ in a punctured
neighborhood of $p$, and $\alpha\cup\beta$ is not geodesic at $p$,
and call these two curves are distinct. As will be seen in the next
section, the angle between two distinct geodesic segments in an
infinite dimensional Teichm\"uller space may still be equal to $0$
or $\pi$, however.

The following result, which follows  directly from Theorem 3.1,
gives a general condition under which there exists the angle between
two geodesic segments in a Teichm\"uller space.

\proclaim{Theorem 4.1} Let $\alpha$ and $\beta$ be two geodesic
segments in $T(R)$ given by the equations $\alpha=\Phi(\mu(t))$ and
$\beta=\Phi(\nu(t))$, $t\in [0, t_0]\, (t_0<1)$, respectively.
Suppose both $\mu(t)$ and $\nu(t)$ are continuous from $[0, t_0]$
into $M(R)$, differentiable at $0$ with $\mu(0)=\nu(0)$, and
$$\tau(\Phi(\mu(0)), \Phi(\mu(t)))=\tau(\Phi(\nu(0)),
\Phi(\nu(t)))=\frac 12\log\frac{1+t}{1-t},\quad t\in [0, t_0].\tag
4.2$$ Then the angle at $\Phi(\mu(0))$ between $\alpha$ and $\beta$
exists, and
$$2\sin\frac{\langle\alpha, \beta\rangle_{\Phi(\mu(0))}}{2}=F(\Phi(\mu(0)),
\Phi'(\mu(0))(\mu'(0)-\nu'(0))).\tag 4.3$$
\endproclaim

We consider a special case of Theorem 4.1. For $\mu_0$, $\mu$ in
$M(R)$, we consider the curve
$$\alpha_{\mu_0, \mu}(t)=\Phi\left(\frac{\delta(\mu_0, \mu)\mu_0(1-\bar\mu_0\mu)+t(\mu-\mu_0)}{\delta(\mu_0, \mu)(1-\bar\mu_0\mu)+t\bar\mu_0(\mu-\mu_0)}\right),\quad t\in [0,
\delta(\mu_0, \mu)],\tag 4.4$$ where $$\delta(\mu_0,
\mu)=\left\|\frac{\mu-\mu_0}{1-\bar\mu_0\mu}\right\|_{\infty}.\tag
4.5$$ We set $\alpha_{\mu}=\alpha_{0, \mu}$ for simplicity.  When
$f^{\mu}\circ (f^{\mu_0})^{-1}$ is extremal, $\alpha_{\mu_0, \mu}$
is a geodesic segment, and
$$\tau([\mu_0], [\alpha_{\mu_0, \mu}(t)])=\frac
12\log\frac{1+t}{1-t}. $$ We call it a standard geodesic segment
joining $[\mu_0]$  to $[\mu]$. By the well-known theorem of
Teichm\"uller, a geodesic segment $\alpha$ beginning at the base
point in a finite dimensional Teichm\"uller space must be a
standard, actually, a Teichm\"uller geodesic segment, that is,
$\alpha=\alpha_{\mu}$ for a so-called Teichm\"uller differential
$\mu=k|\phi|/\phi$ with $0<k<1$, $\phi\in SQ(R)$. While in an
infinite dimensional Teichm\"uller space, a geodesic segment need
not be standard.

The following corollary follows immediately from Theorem 4.1. It
provides an affirmative answer to Problem A posed by Li-Qi [LQ].

\proclaim{Corollary 4.1} Let $\mu_0$, $\mu_1$ and $\mu_2$ be three
 Beltrami coefficients in $M(R)$ such that  $f^{\mu_1}\circ (f^{\mu_0})^{-1}$ and $f^{\mu_2}\circ (f^{\mu_0})^{-1}$ are extremal. Then there exists the
angle at the  point $[\mu_0]$ between the two standard geodesic
segments $\alpha_{\mu_0, \mu_1}$ and $\alpha_{\mu_0, \mu_2}$, and
$$\align &2\sin\frac{\langle\alpha_{\mu_0, \mu_1}, \alpha_{\mu_0, \mu_2}\rangle_{[\mu_0]}}{2}\\
&=F\left(\Phi(\mu_0), \Phi'(\mu_0)
\left(\frac{(\mu_1-\mu_0)(1-|\mu_0|^2)}{\delta(\mu_0,
\mu_1)(1-\bar\mu_0\mu_1)}-\frac{(\mu_2-\mu_0)(1-|\mu_0|^2)}{\delta(\mu_0,
\mu_2)(1-\bar\mu_0\mu_2)}\right)\right).\tag 4.5 \endalign$$ In
particular, when $\mu_0=0$,
$$2\sin\frac{\langle\alpha_{\mu_1}, \alpha_{\mu_2}\rangle_{[0]}}{2}=H\left(\frac{\mu_1}{\|\mu_1\|_{\infty}}-\frac{\mu_2}{\|\mu_2\|_{\infty}}\right).\tag 4.6$$

\endproclaim

In the next section, we will see that  in any infinite dimensional
Teichm\"uller space there do  exist infinitely many  pairs of
geodesic segments (one of which even may be a Teichm\"uller geodesic
segment) between each pair of which the angle does not exist.

 \head 5 An example\endhead

In \S4, we have introduced  the notion of the angle between two
curves in Teichm\"uller spaces and show that such defined angle
exists in a much general situation. A natural question is to
determine whether so-defined angle behaves like that under the
standard Euclidean or hyperbolic geometry. Recall that a geodesic
triangle $\Delta$ in a
 general metric space consists of three distinct geodesic segments, called the sides of $\Delta$, any two of which have precisely one common endpoint. In this section, we will prove the following result.

 \proclaim{Theorem 5.1} Let $R$ be a Riemann surface of
infinite type so that $T(R)$ is infinite dimensional. Given any four
numbers $l$ and $\theta_1$, $\theta_2$, $\theta_3$ with $0<l<\infty$
and $0\le\theta_j\le\pi$ for $j=1, 2, 3$,  there exist infinitely
many geodesic triangles in $T(R)$ each of which has the same three
vertices and a common side and satisfies the property that its three
sides have the same length $l$ while its three angles are equal to
$\theta_1$, $\theta_2$, $\theta_3$ respectively. \endproclaim

Theorem 5.1 implies that the sum of three angles of a geodesic
triangle may be equal to any given number from $0$ to $3\pi$ in an
infinite dimensional Teichm\"uller space. Thus, the geometry of an
infinite dimensional Teichm\"uller space is largely different from
the standard Euclidean or hyperbolic geometry in the view of angle.
This also provides a negative answer to Problem B posed by Li-Qi
[LQ] in the infinite dimensional case. During the proof of Theorem
5.1, we will find out that there do not exist the angles between
infinitely many pairs of geodesic segments (one of which even may be
a Teichm\"uller geodesic segment) in any infinite dimensional
Teichm\"uller space, as stated at the end of $\S$ 4 (see Lemma 5.1
below).

 \vskip 0.3
cm

\noindent {\it Proof of Theorem 5.1}\quad  Let $R$ be a given
Riemann surface of infinite type so that $T(R)$ is infinite
dimensional. Choose an extremal Beltrami coefficient $\mu$ in $M(R)$
which satisfies $|\mu|\equiv k=\frac{e^{2l}-1}{e^{2l}+1}<1$ and
possess a degenerating Hamilton sequence $(\phi_n)$. It is known
there even exist infinitely many Teichm\"uller differentials each of
 which possess a degenerating Hamilton sequence (see [LS]).

Since $(\phi_n)$ is degenerating, we can always choose a sequence of
compact subsets $D_n$ of $R$, and a subsequence of $(\phi_n)$ which
we still denote by $(\phi_n)$, such that
$$D_{n-1}\subset D_n,\quad R=\cup_{n=1}^{\infty}D_n,\tag 5.1$$
and
$$\iint_{D_n\setminus D_{n-1}}|\phi_n|dxdy=1+o(1),\quad
n\to\infty.\tag 5.2$$

We need to consider the inverse map $(f^{\mu})^{-1}$, and denote by
$\mu^*$ its Beltrami coefficient. Since $\mu$ is extremal with a
degenerating Hamilton sequence, $\mu^*$ is also extremal, and has a
degenerating Hamilton sequence, which we denote by $(\phi^*_n)$. Set
$D^*_n=f^{\mu}(D_n)$. Then $(D^*_n)$ is a sequence of compact
subsets of $R^{\mu}$. We now choose subsequences of $(\phi^*_n)$ and
$(D_n^*)$, which we still denote by $(\phi^*_n)$ and $(D_n^*)$, such
that
$$D^*_{n-1}\subset D^*_n,\quad R^{\mu}=\cup_{n=1}^{\infty}D^*_n,\tag 5.3$$
and
$$\iint_{D^*_n\setminus D^*_{n-1}}|\phi^*_n|dxdy=1+o(1),\quad
n\to\infty.\tag 5.4$$

We list two basic properties of these constructions. Set
$R_1=\cup_{n=0}^{\infty}(D_{2n+1}\setminus D_{2n})$
$(D_0=\varnothing)$, $R_2=R\setminus R_1$,
$R_1^{\mu}=\cup_{n=0}^{\infty}(D^*_{2n+1}\setminus D^*_{2n})$
$(D^*_0=\varnothing)$, $R^{\mu}_2=R^{\mu}\setminus R^{\mu}_1$. Let
$\chi$ denote the characteristic function of a set. Then for any two
numbers $c_1$ and $c_2$, we have the following two statements:
$$H_R((c_1\chi_{R_1}+c_2\chi_{R_2})\mu)=\|(c_1\chi_{R_1}+c_2\chi_{R_2})\mu\|_{\infty}=k\max(|c_1|,
|c_2|).\tag $P_1$
$$
$$H_{R^{\mu}}((c_1\chi_{R^{\mu}_1}+c_2\chi_{R^{\mu}_2})\mu^*)=\|(c_1\chi_{R^{\mu}_1}+c_2\chi_{R^{\mu}_2})\mu^*\|_{\infty}=k\max(|c_1|,
|c_2|).\tag $P_2$
$$

\noindent {\bf Step 1}\quad Constructing a geodesic segment between
$\Phi(\mu)$ and $\Phi(\chi_{R_1}\mu)$ \vskip 0.3 cm For simplicity,
set $\mu_1=\chi_{R_1}\mu$. $(P_1)$ implies that $\mu_1$ is extremal,
and $\|\mu_1\|_{\infty}=k$. To construct a geodesic segment between
$[\mu]$ and $[\mu_1]$, we adapt some discussion from Li [Li3] and
the  second-named author [Sh]. Define $\mu_t$ in $M(R)$ joining $\mu$
to $\mu_1$ as follows:
$$\mu_t=\left(\sigma(t)\chi_{R_1}+
\frac{k-t}{k(1-kt)}\chi_{R_2}\right)\mu,\tag 5.5
$$
where $\sigma(t)$ is a continuous function of $t$ in $[0, k]$
satisfying the following condition:
$$\max\left\{\frac{k-t}{k(1-kt)},
\frac{t}{k}\right\}\le\sigma(t)\le\min\left\{\frac{k+t}{k(1+kt)},
\frac{2k-(1+k^2)t}{k(1+k^2-2kt)}\right\}.\tag 5.6$$ Clearly,
$\|\mu_t\|_{\infty}=\sigma(t)k$. Using $(P_1)$ again, we see that
each $\mu_t$ is extremal. We first prove that  $\Phi(\mu_t)$, $t\in [0, k]$, is a
geodesic segment between $[\mu]$ and $[\mu_1]$.

In fact, if we set $f_t=f^{\mu_t}$, then the complex dilatation
$\tilde\mu_t$ of $f_t\circ f_0^{-1}$ is
$$\tilde\mu_t=\left(\frac{\mu_t-\mu}{1-\bar\mu\mu_t}\cdot\frac{{\partial}_z
f_0}{\overline{{\partial}_z f_0}}\right)\circ f_0^{-1}=\left(
 \frac{1-\sigma(t)}{1-k^2\sigma(t)}\chi_{R^{\mu}_1}+
\frac{t}{k}\chi_{R^{\mu}_2}\right)\mu^*.\tag 5.7
$$
A direct but tedious computation from (5.6) yields
$\|\tilde\mu_t\|_{\infty}=t.$ On the other hand, by $(P_2)$ we
conclude that $\tilde\mu_t$ is extremal. So we get
$$\tau([\mu], [\mu_t])=\frac 12\log\frac{1+t}{1-t}.\tag 5.8$$

Now the Beltrami coefficient $\tilde\nu_t$ of $f_t\circ f_1^{-1}$ is
$$\tilde\nu_t=\left(\frac{\mu_t-\mu_1}{1-\bar\mu_1\mu_t}\cdot\frac{{\partial}_z
f_1}{\overline{{\partial}_z f_1}}\right)\circ f_1^{-1}.\tag 5.9$$ By the
definition of $\mu_t$  and the inequality (5.6) we get
$$\|\tilde\nu_t\|_{\infty}=\frac{k-t}{1-kt}.\tag 5.10$$
On the other hand, by (5.8) and (5.10),  $$\tau([\mu_1],
[\mu_t])\ge\tau([\mu], [\mu_1])-\tau([\mu], [\mu_t])=\frac
12\log\left(\frac{1+k}{1-k}\cdot\frac{1-t}{1+t}\right)=\frac
12\log\frac{1+\|\tilde\nu_t\|_{\infty}}{1-\|\tilde\nu_t\|_{\infty}}.$$
Therefore, $\tilde\nu_t$ is extremal, and
$$\tau([\mu_1],
[\mu_t])=\frac
12\log\frac{1+\|\tilde\nu_t\|_{\infty}}{1-\|\tilde\nu_t\|_{\infty}}.\tag
5.11$$ Consequently, $$\tau([\mu], [\mu_t])+\tau([\mu_1],
[\mu_t])=\tau([\mu], [\mu_1]),\quad t\in [0, k].$$ This implies that
$\Phi(\mu_t)$, $t\in [0, k]$, is a geodesic segment between $[\mu]$
and $[\mu_1]$.

Now we consider the standard  geodesic segments
$\alpha_{\mu}=\Phi(t/k\mu)$, $\alpha_{\mu_1}=\Phi(t/k\mu_1)$ and the
above-constructed geodesic segment $\beta_{\sigma}=\Phi(\mu_t)$,
$t\in [0, k]$. Then they have the same length $\frac
12\log\frac{1+k}{1-k}=l.$ Clearly,
$\alpha_{\mu}\cap\beta_{\sigma}=\{[\mu]\}$,
$\alpha_{\mu_1}\cap\beta_{\sigma}=\{[\mu_1]\}$. Before we discuss
the existence of the angles at $[\mu]$ and $[\mu_1]$, we point out when
$\alpha_{\mu}\cup\beta_{\sigma}$ is (not) geodesic at $[\mu]$, and
$\alpha_{\mu_1}\cup\beta_{\sigma}$ is (not) geodesic at $[\mu_1]$.
Since both $\alpha_{\mu}$ and $\beta_{\sigma}$ are geodesic
segments, $\alpha_{\mu}\cup\beta_{\sigma}$ is not geodesic at
$[\mu]$ if and only if, as $t\to 0_+$,
$$\tau([0], [\mu])+\tau([\mu], [\mu_t])>\tau([0], [\mu_t]),$$
which implies by (5.8) that
$$\sigma(t)<\frac{k+t}{k(1+kt)},\quad t\to 0_+.\tag 5.12$$
Similarly, $\alpha_{\mu_1}\cup\beta_{\sigma}$ is not geodesic at
$[\mu_1]$ if and only if
$$\sigma(t)<\frac{2k-(1+k^2)t}{k(1+k^2-2kt)},\quad t\to k_-.\tag
5.13$$ Under these two conditions, $\alpha_{\mu}$, $\alpha_{\mu_1}$
and $\beta_{\sigma}$ are distinct geodesic segments. In the
following, we assume that $\sigma$ satisfies (5.6), (5.12) and
(5.13). Corollary (4.1) and $(P_1)$ imply that
$$\langle\alpha_{\mu},
\alpha_{\mu_1}\rangle_{[0]}=2\arcsin\frac{H(\mu-\mu_1)}{2k}=2\arcsin\frac{H(\chi_{R_2}\mu)}{2k}=\frac{\pi}{3}.$$

We end this step by pointing out that $\sigma(t)\equiv 1$ meets all
the conditions (5.6), (5.12) and (5.13). In this case
$\beta_{\sigma}$ is the standard geodesic segment $\alpha_{\mu,
\mu_1}$. This will be essential in our final step to construct the
desired geodesic triangle. We also point out that there are
infinitely many continuous functions $\sigma$ satisfying the
conditions (5.6), (5.12) and (5.13), and two different such
functions determine two different geodesic segments between $[\mu]$
and $[\mu_1]$.

\vskip 0.3 cm \noindent {\bf Step 2}\quad On the existence of the
angle at $[\mu]$ between $\alpha_{\mu}$ and $\beta_{\sigma}$

\proclaim{Lemma 5.1} There exists  the angle at $[\mu]$ between
$\alpha_{\mu}$ and $\beta_{\sigma}$ if and only if $\sigma$ is
differentiable at $0$. Furthermore, we can choose $\sigma$ so that
the angle at $[\mu]$ between $\alpha_{\mu}$ and $\beta_{\sigma}$
attains any given number from $0$ to $\pi$.\endproclaim \demo{Proof}
We first assume $\sigma$ is differentiable at $0$. (5.6) implies
that $|\sigma'(0)|\le (1-k^2)/k$. Then $\mu_t$ is differentiable
with $\mu_0=\mu$, and
$$\mu'(0)=\lim_{t\to 0_+}\frac{\mu_t-\mu_0}{t}=\left(
 \sigma'(0)\chi_{R_1}+
\frac{k^2-1}{k}\chi_{R_2}\right)\mu.\tag 5.14
$$
Interchanging the endpoints, $\alpha_{\mu}=\Phi(\nu_t)$, with
$$\nu_t=\frac{k-t}{k(1-kt)}\mu,\quad t\in [0, k].$$ Then $\nu_t$ is
differentiable with $\nu_0=\mu$, and
$$\nu'(0)=\lim_{t\to 0_+}\frac{\nu_t-\mu}{t}=\frac{k^2-1}{k}\mu.\tag 5.15$$
It is easy to see that
$$\tau([\mu], [\nu_t])=\frac 12\log\frac{1+t}{1-t}.\tag 5.16$$ By
(5.8), (5.14-16), we find out that $\mu_t$ and $\nu_t$ satisfy the
assumption in Theorem 4.1. Consequently, there exists the angle at
$[\mu]$ between $\alpha_{\mu}$ and $\beta_{\sigma}$, and
$$2\sin\frac{\langle\alpha_{\mu},
\beta_{\sigma}\rangle_{[\mu]}}{2}=F(\Phi(\mu),
\Phi'(\mu)(\mu'(0)-\nu'(0))).$$ Noting that
$$\mu'(0)-\nu'(0)=\left(\frac{1-k^2}{k}+\sigma'(0)\right)\chi_{R_1}\mu,$$
we obtain
$$\align
F(\Phi(\mu),
\Phi'(\mu)(\mu'(0)-\nu'(0)))&=\left(\frac{1-k^2}{k}+\sigma'(0)\right)F(\Phi(\mu),
\Phi'(\mu)(\chi_{R_1}\mu))\\
&=\left(\frac{1}{k}+\frac{\sigma'(0)}{1-k^2}\right)H_{R^{\mu}}\left(\left(\chi_{R_1}\mu\cdot\frac{{\partial}_z
f_0}{\overline{{\partial}_z f_0}}\right)\circ f_0^{-1}\right)\\
&=\left(\frac{1}{k}+\frac{\sigma'(0)}{1-k^2}\right)H_{R^{\mu}}(\chi_{R^{\mu}_1}\mu^*)=1+\frac{k\sigma'(0)}{1-k^2}.\endalign
$$
Consequently, $$ 2\sin\frac{\langle\alpha_{\mu},
\beta_{\sigma}\rangle_{[\mu]}}{2}=1+\frac{k\sigma'(0)}{1-k^2}.\tag
5.17$$
 Since $|\sigma'(0)|\le (1-k^2)/k$, we see that
$\langle\alpha_{\mu}, \beta_{\sigma}\rangle_{[\mu]}\in [0, \pi]$.

To show that $\langle\alpha_{\mu}, \beta_{\sigma}\rangle_{[\mu]}$
may attain any number from $0$ to $\pi$, it is sufficient to show
that for any number $\delta$ with $|\delta|\le (1-k^2)/k$, there
exists a continuous function $\sigma_{\delta}(t)$ in $[0, k]$ which
satisfies the inequalities (5.6), (5.12) and is differentiable at
zero with $\sigma_{\delta}'(0)=\delta$. Actually, by (5.12), we only
need to find $\sigma_{\delta}(t)$ when $t\to 0_+$. When
$\delta=(1-k^2)/k$, we choose
$$\sigma_{\delta}(t)=1+\frac{1-k^2}{k}t-(1-k^2)t^2,\quad t\to 0_+,$$
when $-(1-k^2)/k\le\delta<(1-k^2)/k$, we choose
$$\sigma_{\delta}(t)=1+\delta t,\quad t\to 0_+.$$

Conversely, suppose that there exists the angle at $[\mu]$ between
$\alpha_{\mu}$ and $\beta_{\sigma}$. By $(5.8)$ and $(5.16)$ we
conclude that $$\lim_{t\to 0_+}\frac{\tau([\mu_t],
[\nu_t])}{t}=2\sin\frac{\langle\alpha_{\mu},
\beta_{\sigma}\rangle_{[\mu]}}{2}.$$ Let $\mu^*_t$ denote the
Beltrami coefficient of $f^{\mu_t}\circ(f^{\nu_t})^{-1}$. It is
routine to show that $\mu^*_t$ is extremal, and
$$k_0(\mu^*_t)=\left\|\frac{k\sigma(t)-\frac{k-t}{1-kt}}{1-\frac{k-t}{1-kt}k\sigma(t)}\right\|_{\infty}.\tag 5.18$$
Then, $$\tau([\mu_t], [\nu_t])=\frac
12\log\frac{1+k_0(\mu^*_t)}{1-k_0(\mu^*_t)}=k_0(\mu^*_t)+o(t)=\frac{(\sigma(t)-1)k}{1-k^2}+t+o(t),\quad
t\to 0_+.$$ Consequently, $\sigma$ is differentiable at $0$, and
$$\sigma'(0)=\frac{1-k^2}{k}\left(2\sin\frac{\langle\alpha_{\mu},
\beta_{\sigma}\rangle_{[\mu]}}{2}-1\right).$$

Finally, we need to find a continuous function $\sigma$ which
satisfies (5.6), (5.12) but is not differentiable at $0$ so that the
angle at $[\mu]$ between $\alpha_{\mu}$ and $\beta_{\sigma}$ does
not exist. In fact,  the following function works:
$$\sigma(t)=1+\frac{1-k^2}{2k}t\sin^2\frac 1t,\quad t\to 0_+.$$

\enddemo

\vskip 0.3 cm \noindent {\bf Step 3}\quad On the existence of the
angle at $[\mu_1]$ between $\alpha_{\mu_1}$ and $\beta_{\sigma}$

\proclaim{Lemma 5.2} There exists  the angle at $[\mu_1]$ between
$\alpha_{\mu_1}$ and $\beta_{\sigma}$ if and only if $\sigma$ is
differentiable at $k$. Furthermore, we can choose $\sigma$ so that
the angle at $[\mu_1]$ between $\alpha_{\mu_1}$ and $\beta_{\sigma}$
attains any given number from $0$ to $\pi$.\endproclaim \demo{Proof}
By the same reasoning as in Step 2, we can prove that the angle at
$[\mu_1]$ between $\alpha_{\mu_1}$ and $\beta_{\sigma}$ exists if
and only if $\sigma$ is differentiable at $k$.  In this case,
$|\sigma'(k)|\le 1/k$,  and
$$2\sin\frac{\langle\alpha_{\mu_1},
\beta_{\sigma}\rangle_{[\mu_1]}}{2}=1-k\sigma'(k).\tag 5.19$$ We
omit the details here. Consequently, $\langle\alpha_{\mu_1},
\beta_{\sigma}\rangle_{[\mu_1]}\in [0, \pi]$, and it attains  any
number from $0$ to $\pi$. As above, it is sufficient to show that
for any number $\eta$ with $|\eta|\le 1/k$, there exists a
continuous function $\sigma_{\eta}(t)$ in $[0, k]$ which satisfies
the inequalities (5.8), (5.13) and is differentiable at $k$ with
$\sigma_{\eta}'(k)=\eta$. This time, by (5.13), we only need to find
$\sigma_{\delta}(t)$ when $t\to k_-$.When $\eta=-1/k$, we choose
$$\sigma_{\eta}(t)=1-\frac{t-k}{k}-\frac{2(t-k)^2}{1-k^2},\quad t\to k_-,$$
when $-1/k<\eta\le 1/k$, we choose
$$\sigma_{\eta}(t)=1+\eta (t-k),\quad t\to k_-.$$
\enddemo

\noindent {\bf Step 4}\quad Constructing the  triangle \vskip 0.3 cm

We have proved that
$\alpha_{\mu}\cup\alpha_{\mu_1}\cup\beta_{\sigma}$ is a geodesic
triangle such that its two  angles at $[\mu]$ and $[\mu_1]$ can be
equal to any two given numbers from $0$ to $\pi$ for an appropriate
function $\sigma$, but the third angle at $[0]$ is fixed and equal
to $\pi/3$. We now show that $\alpha_{\mu_1}$ can be modified in a
neighborhood of $[0]$ so that the new angle at $[0]$ can be equal to
any given number from $0$ to $\pi$.

Recall that the existence and the value of the angle between two
geodesic segments are preserved under an allowable map. By (5.7) we
see that $\tilde\mu_1=\chi_{R^{\mu}_2}\mu^*$. Replacing $\mu$,
$\mu_1$ and $R_1$ by $\mu^*$, $\tilde\mu_1$ and $R^{\mu}_2$
respectively in the three steps above, we obtain a geodesic segment
$\gamma_{\tilde\sigma}$ in $T(R^{\mu})$ joining $[\mu^*]$ to
$[\tilde\mu_1]$ such that the angle at $[\mu^*]$ between
$\gamma_{\tilde\sigma}$ and $\alpha_{\mu^*}$ exists and equals any
given number from $0$ to $\pi$, and $\gamma_{\tilde\sigma}$
coincides with the standard geodesic segment $\alpha_{\mu^*,
\tilde\mu_1}$ in a neighborhood of $[\tilde\mu_1]$ (as remarked at
the end of Step 1). By the allowable map $A_{\mu}: T(R)\to
T(R^{\mu})$, $\gamma_{\tilde\sigma}$ becomes a geodesic segment
$\tilde\gamma_{\tilde\sigma}$ in $T(R)$ joining $[0]$ to $[\mu_1]$
such that the angle at $[0]$ between $\tilde\gamma_{\tilde\sigma}$
and $\alpha_{\mu}$ exists and equals any given number from $0$ to
$\pi$, and $\tilde\gamma_{\tilde\sigma}$ coincides with the standard
geodesic segment $\alpha_{\mu_1}$ in a neighborhood of $[\mu_1]$.
Now let $\Delta_0$ denote the geodesic triangle
$\alpha_{\mu}\cup\beta_{\sigma}\cup\tilde\gamma_{\tilde\sigma}$ with
vertices $[0]$, $[\mu]$ and $[\mu_1]$. Then the angles at these
three vertices exist and equal any three given numbers from $0$ to
$\pi$. This finishes our construction and completes the proof of
Theorem 5.1.

\vskip 0.3 cm \noindent {\bf Remark 1} In a finite dimensional
Teichm\"uller space, the angle between two distinct geodesic
segments exists and is always positive and less than $\pi$. It is
not known whether the sum of three angles of a geodesic triangle in
a  finite dimensional Teichm\"uller space is less than $\pi$. This
seems to be a difficult problem.

\vskip 0.3 cm \noindent {\bf Remark 2} Masur [Ma] proved that any
 Teichm\"uller space of finite dimension $(\ge  2)$ does not have  negative curvature.  In general, a metric space
$(X, d)$ is said to have negative curvature if for any geodesic
triangle $\Delta$ in $X$ with vertices $A$, $B$ and $C$, $d(B,
C)>2d(\tilde B, \tilde C)$, where $\tilde B$ is the midpoint of the
side $\overline{AB}$ between $A$ and $B$,  and $\tilde C$ is the
midpoint of the side $\overline{AC}$ between $A$ and $C$. In fact,
by Masur's discussion,  for any
 Teichm\"uller space $T(R)$ of finite dimension $(\ge  2)$ and any
 number $\epsilon>1$, there exists a geodesic triangle $\Delta$ in $T(R)$ with vertices $A$, $B$ and $C$ such that $d(B, C)\le\epsilon d(\tilde
B, \tilde C)$. It is of interest to determine whether one can take
$\epsilon=1$ here. Our proof of Theorem 5.1 shows that this is the
case when $T(R)$ is infinite dimensional. Thus, an infinite
dimensional Teichm\"uller space is far from having negative
curvature.

In fact, we may consider the triangle
$\Delta=\alpha_{\mu}\cup\alpha_{\mu_1}\cup\beta_{\sigma}$. Here we
assume that
$$\sigma(t_0)=\frac{2+\sqrt{1-k^2}}{1+k^2+\sqrt{1-k^2}}\tag 5.20$$
with
$$t_0=\frac{k}{1+\sqrt{1-k^2}}.\tag 5.21$$
By (5.8), (5.16) and (5.20), we obtain
$$\tau([\mu], [\mu_{t_0}])=\tau([\mu], [\nu_{t_0}])=\frac
14\log\frac{1+k}{1-k}=\frac l2.$$ Thus, $[\mu_{t_0}]$ and
$[\nu_{t_0}]$ are the midpoints of $\beta_{\sigma}$ and
$\alpha_{\mu}$, respectively. Now it follows from (5.18), (5.20) and
(5.21) that
$$\tau([\mu_{t_0}], [\nu_{t_0}])=\frac
12\log\frac{1+k_0(\mu^*_{t_0})}{1-k_0(\mu^*_{t_0})}=\frac
12\log\frac{1+k}{1-k}=l=\tau([\mu_1], [0]).$$

\vskip 0.2 cm

\noindent {\bf Acknowledgements} \quad
 The authors would like to thank the  referee for useful advice.

 
\enddocument